\documentclass[titlepage,11pt]{article}
\usepackage{latexsym, tikz}
\usepackage{amsfonts}
\usepackage{amsmath}
\usepackage{textcomp}
\usetikzlibrary{graphs}
\usepackage{float}

\newcommand\blackslug{\hbox{\hskip 1pt \vrule width 4pt height 8pt depth 1.5pt
        \hskip 1pt}}
\newcommand\bbox{\hfill \quad \blackslug \bigbreak}
\def\dd{\hbox{-}}
\def\cc{\hbox{-}\cdots\hbox{-}}
\def\ll{,\ldots,}

\title{Detecting an odd hole}
\author{Maria Chudnovsky\thanks{This material is based upon work supported in part by the U. S. Army
Research Office under grant   number W911NF-16-1-0404, and supported by  NSF grant DMS 1763817.}\\
Princeton University, Princeton, NJ 08544
\\
\\
Alex Scott\thanks{Supported by a Leverhulme Trust Research
Fellowship.}\\
Mathematical Institute, University of Oxford, Oxford OX2 6GG, UK
\\
\\
Paul Seymour\thanks{Supported by NSF grant DMS-1800053.}\\
Princeton University, Princeton, NJ 08544
\\
\\
Sophie Spirkl\thanks{This material is based upon work supported by the National Science
Foundation under Award No. DMS-1802201.}\\
Rutgers University, New Brunswick, NJ 08901}

\date{December 21, 2018; revised \today}

\newtheorem{thm}{}[section]

\newcommand{\Proof}{\noindent{\bf Proof.}\ \ }

\begin{document}
\maketitle
\begin{abstract}
We give a polynomial-time algorithm to test whether a graph contains an induced cycle with length more than three and odd.
\end{abstract}

\section{Introduction}
All graphs in this paper are finite and have no loops or parallel edges. A {\em hole} of $G$ is an induced subgraph of $G$
that is a cycle
of length at least four. In this paper we give an algorithm to test whether a graph $G$ has an odd hole, with running time polynomial 
in $|G|$. ($|G|$ denotes the number of vertices of $G$.)

The study of holes, and particularly odd holes, grew from Claude Berge's ``strong perfect graph conjecture''~\cite{Berge}, that if a graph and its
complement both have no odd holes, then its chromatic number equals its clique number. For many years, this was an open question, as was the question
of finding a polynomial-time algorithm to test if a graph is perfect. (``Perfect'' means every induced subgraph has chromatic number
equal to clique number.) Both questions were settled at about the same time: in the early 2000's,
two of us, with Robertson and 
Thomas~\cite{SPGT}, proved Berge's conjecture, and also, with Cornu\'ejols, Liu and Vu\v{s}kovi\'{c}~\cite{bergealg}, gave a polynomial-time algorithm to test if a graph
has an odd hole or odd antihole, thereby testing perfection.
(An {\em antihole} of $G$ is an induced subgraph whose complement graph is a hole in the complement graph of $G$.)

Excluding both odd holes and odd antiholes is a combination that works well, and has ``deep'' structural consequences. Just excluding
odd holes loses this advantage, and graphs with no odd holes seem to be much less well-structured than
perfect graphs, in some vague sense.  It was only recently that two of us~\cite{oddholes} proved
that if a graph has no odd holes then its chromatic number is bounded by a
function of its clique number, resolving an old conjecture of Gy\'arf\'as~\cite{gyarfas}.
For a long time, the question of testing for odd holes has
remained open: while we could test for the presence of a odd hole or
antihole in polynomial time, we were unable to separate the test for odd holes from the test
for odd antiholes, and the complexity of testing for an odd hole remained
open. Indeed, it seemed quite likely that testing for an odd hole was NP-complete; for instance,
D. Bienstock~\cite{bienstock,bienstock2} showed that testing if a graph has an odd hole containing a given vertex is NP-complete.
But in fact we can test for odd holes in polynomial time.

The main result of this paper is the following:

\begin{thm}\label{mainthm}
There is an algorithm with the following specifications:
\begin{description}
\item [Input:] A graph $G$.
\item [Output:] Determines whether $G$ has an odd hole.
\item [Running time:] $O(|G|^{9})$.
\end{description}
\end{thm}
Remarkably, the algorithm (or rather, the proof of its correctness) is considerably simpler than the algorithm of~\cite{bergealg}, 
and so not only
solves an open question, but gives a better way to test if a graph is perfect. Its running time is the same as the old algorithm 
for testing perfection.

Like the algorithm of~\cite{bergealg} (and several other
algorithms that test for the presence of an induced subdivision of a fixed
subgraph) the new algorithm has three stages.
Say we are looking for an induced subgraph of ``type T'' in the input graph $G$. 
\begin{itemize}
\item The algorithm first tests for certain
``easily-detected configurations''. 
These are configurations that can be efficiently detected and whose presence guarantees that $G$ contains an induced subgraph of type T.
\item The third step is an algorithm that tries to find a subgraph of type T directly; it would not be expected to work on a general input graph,
but it would detect a subgraph of type T if there happens to be a particularly ``nice'' one in the input graph. For instance, in \cite{bergealg} 
we were looking for an odd hole, and the method was, try all triples of vertices and join them by three shortest paths, and see if they form an odd hole.
This would normally not work, but it would work if there was an odd hole that was both the shortest odd hole in the graph and there were no vertices 
in the remainder
of the graph with many neighbours in this hole. 
\item The second step is to prepare the
input for the third step; this step is called ``cleaning'', and is where the main challenges lie.
In the cleaning step, the algorithm generates a ``cleaning list'', polynomially many subsets
$X_1, \ldots, X_k$ of the vertex
set of the input graph $G$, with the property
that if $G$
does in fact contain an induced subgraph of type T, then 
for some $i \in \{1, \ldots, k\}$, such a subgraph can be
found in $G \setminus X_i$  using the
method of step 3. This usually means that if $G$ contains an induced subgraph of type T, then there is one, say $H$, such that
some $X_i$ contains
all the vertices of $G\setminus V(H)$ that have many neighbours in $H$, and
deleting $X_i$ leaves $H$ intact.
\end{itemize}

Cleaning was first used  by Conforti and Rao~\cite{conforti} to recognize
linear balanced matrices, and subsequently by Conforti,
Cornu\'ejols, Kapoor and Vu\v{s}kovi\'{c}~\cite{kapoor} to test for even
holes, as well as in \cite{bergealg}.
It then became a standard tool in induced subgraph detection
algorithms~\cite{prismtheta, net,kawa}.
This is the natural approach to try to test for an odd hole, and it seemed to have been explored thoroughly;
but not thoroughly enough, as we shall see. We have found a novel method of cleaning that works remarkably well, and
and will have further applications. 

In both the old algorithm to check perfection,  and the new algorithm of this paper,
we first test whether $G$ contains a ``pyramid'' or ``jewel'' (easily-detected 
induced subgraphs that would imply the presence of an odd hole) and we may assume that it does not. 
Then we try to search for
an odd hole directly, and to do so we exploit the properties of an odd hole of minimum length, a so-called {\em shortest odd hole}. 
(These properties hold in graphs with no pyramid or jewel, but not in general graphs, which is why we test for pyramids and jewels first.)

Let $C$ be a shortest odd hole.
A vertex $v\in V(G)$ is {\em $C$-major} if there is no three-vertex path of $C$ containing all the neighbours of $v$ in $V(C)$
(and consequently $v\notin V(C)$); and $C$ is {\em clean} (in $G$) if no vertices of $G$ are $C$-major. 
If $G$ has a shortest odd hole 
that is clean, then it is easy to detect that $G$ has an odd hole, and this was done in theorem 4.2 of~\cite{bergealg}.
More exactly, there is a poly-time algorithm that either finds an odd hole, or deduces that no shortest odd hole of $G$ is clean.
Call this {\em procedure P}.
The complicated part of the algorithm in~\cite{bergealg} was generating a cleaning list,
a list of polynomially-many subsets of $V(G)$,
such that if $G$ has an odd hole, then there is a shortest odd hole $C$ and some $X$ in the list such that
$X\cap V(C)=\emptyset$ and every $C$-major vertex belong to $X$. Given that, we take the list $X_1\ll X_k$ say,
and for each of the polynomially-many graphs $G\setminus X_i$, we run procedure P on it. If it ever finds an odd hole, then $G$
has an odd hole and we are done. If not, then $G$ has no odd hole and again we are done.

So, the key is generating the cleaning list. For this, \cite{bergealg} uses
\begin{thm}\label{completedge}
For every graph $G$ not containing any ``easily-detected configuration'', if $C$ is a shortest odd hole in $G$, and $X$
is an anticonnected set of $C$-major vertices, then there is an edge $uv$
of $C$ such that $u,v$ are both $X$-complete. 
\end{thm}
({\em Anticonnected}  means connected in the complement, and {\em $X$-complete} means adjacent to every vertex in $X$.)
But to arrange that \ref{completedge} is true,
it is necessary to expand the definition of ``easily-detected configuration'' to include some new configurations.
It remains true that if one is present then the graph has an odd hole or odd antihole and we can stop,
and if they are not present then \ref{completedge} is true. The problem is, 
if one of these new
easily-detected configurations is present, it guarantees that $G$ contains an odd hole or an odd antihole,
but not necessarily an odd hole.

But there is a simpler way.
Here is a rough sketch of a new procedure to 
clean a shortest odd hole $C$ in a graph $G$ with no pyramid or jewel. Let $x$ be a $C$-major vertex such that there is a gap in $C$ between two neighbours 
of $x$, as long as possible. (We can assume there is one.)
Let the neighbours of $x$ at the ends of this gap be $d_1,d_2$; thus there is a path $D$ of $C$ between
$d_1,d_2$ such that every $C$-major vertex either has a neighbour in its interior, or is adjacent to both ends. We can assume that the $C$-distance
between $d_1,d_2$ is at least three. Also, 
we have a theorem that there is an edge $f$ of $C$ such that every $C$-major vertex nonadjacent to $x$ is adjacent to one of the ends of $f$.

For the algorithm, what we do is: we guess $x,d_1,d_2$ and $f$ (more precisely, we enumerate all possibilities for them). Eventually we will guess 
correctly. We also guess the two vertices neighbouring $f$ in $C$, say $c_1,c_4$, where $f=c_2c_3$ and $c_1,c_2,c_3,c_4$ are in order in $C$. 
When we guess correctly, every $C$-major vertex either
\begin{itemize}
\item is adjacent to both $d_1,d_2$; or 
\item is different from $c_1,c_2,c_3,c_4$ and is adjacent to one of $c_2,c_3$; or
\item is adjacent to $x$ and has a neighbour in the interior of $D$.
\end{itemize}
We can safely delete all common neighbours of $d_1,d_2$ except $x$; deleting these vertices will not remove any vertices
of $C$. Also, we can safely delete all vertices different from $c_1,c_2,c_3,c_4$  that are adjacent to one of $c_2,c_3$.
So now in the graph that remains after these deletions, say $G'$,
all $C$-major vertices different from $x$ satisfy the third bullet above. 

We do not know the path $D$, and so we cannot immediately identify 
the set of vertices satisfying
the third bullet. (For this sketch, let us assume that $D$ has length less than half that of $C$; if it is longer there is a slight complication.)
But we know (it is a theorem of \cite{bergealg}) that $D$ is a shortest path between $d_1,d_2$ in the graph obtained from $G'$ by deleting all $C$-major 
vertices; and so it is also a shortest path between $d_1,d_2$ in the graph $G''$ obtained from $G'$ by deleting $x$ and all its neighbours (except $d_1,d_2$).
The algorithm computes $G''$, and then 
finds the union of the interiors of the vertex sets of all shortest paths between $d_1,d_2$ in $G''$, say $F$. 
It is another theorem of~\cite{bergealg}
that no vertex of $C\setminus V(D)$ has a neighbour in $F$; so it is safe to delete from $G'$ all vertices of $G'$  except $d_1,d_2$ that are not in $F$ 
and have a neighbour in $F$. But then we have deleted all the $C$-major vertices, and now we just test for a clean shortest odd hole.

In an earlier version of this paper, we proved the result by a more complicated method that also seems to us novel and worth recording.
It was necessary to first test for two more easily-detected configurations; but then, instead
of constructing the set $F$ above, the algorithm just guesses the component ($F'$ say) of $G''$ that contains the interior of $D$, and deletes
all neighbours of $x$ that have neighbours in this component except $d_1,d_2$. This might delete some of the hole $C$, but we proved a theorem that enough of $C$
remains that we can still use it in an algorithm to detect an odd hole. In particular, there is an odd path $P$ of $C$ of length at least three,
with both ends adjacent to $x$, such that the ends of $P$ both have neighbours in $F'$ and its internal vertices do not; and we can exploit this to detect
the presence of an odd hole.

\section{The easily-detected configurations}

Let $v_0\in V(G)$, and for $i = 1,2,3$ let $P_i$ be an induced path of $G$ between $v_0$ and $v_i$, such that

\begin{itemize}
\item $P_1,P_2,P_3$ are pairwise vertex-disjoint except for $v_0$;
\item $v_1,v_2,v_3\ne v_0$, and at least two of $P_1,P_2,P_3$ have length at least two;
\item $v_1,v_2,v_3$ are pairwise adjacent; and
\item for $1\le i<j\le 3$, the only edge between $V(P_i)\setminus \{v_0\}$ and $V(P_j)\setminus \{v_0\}$ is the edge $v_iv_j$.
\end{itemize}
\begin{figure}[H]
\centering

\begin{tikzpicture}[scale=0.8,auto=left]
\tikzstyle{every node}=[inner sep=1.5pt, fill=black,circle,draw]

\node (v0) at (0,0) {};
\node (v1) at (-1,-3.5) {};
\node (v2) at (0,-3) {};
\node (v3) at (1,-3.5) {};
\tikzstyle{every node}=[]
\draw (v1) node [below]           {\footnotesize$v_1$};
\draw (v2) node [below]           {\footnotesize$v_2$};
\draw (v3) node [below]           {\footnotesize$v_3$};
\draw (v0) node [above]           {\footnotesize$v_0$};

\draw (v1) -- (v2);
\draw (v1) -- (v3);
\draw (v2) -- (v3);
\draw[dashed] (v0) -- (v1);
\draw[dashed] (v0) -- (v2);
\draw[dashed] (v0) -- (v3);

\end{tikzpicture}

\caption{A pyramid. Throughout, dashed lines represent paths, of indeterminate length.} \label{fig:pyramid}
\end{figure}
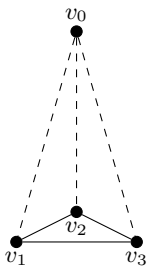

We call the subgraph induced on $V(P_1\cup P_2\cup P_3)$ a {\em pyramid}, with {\em apex} $v_0$ and {\em base} $\{v_1,v_2,v_3\}$.
If $G$ has a pyramid then $G$ has an odd hole (because two of the paths $P_1,P_2,P_3$ have the same length modulo two, and they induce
an odd hole). It is shown in theorem 2.2 of~\cite{bergealg} that:

\begin{thm}\label{testpyramid}
There is an algorithm with the following specifications:
\begin{description}
\item [Input:] A graph $G$.
\item [Output:] Determines whether there is a pyramid in $G$.
\item [Running time:] $O(|G|^{9})$.
\end{description}
\end{thm}

If $X\subseteq V(G)$, we denote the subgraph of $G$ induced on $X$ by $G[X]$. If $X$ is a vertex or edge of $G$, or a set of vertices 
or a set of edges of $G$, we denote by $G\setminus X$ the graph obtained from $G$ by deleting $X$.
Thus, for instance, if $b_1b_2$ is an edge of a hole $C$, then $C\setminus \{b_1,b_2\}$ and $C\setminus b_1b_2$ are both paths,
but one contains $b_1,b_2$ and the other does not. If $P$ is a path, 
the {\em interior} of $P$ is the set of vertices of the path $P$ that are not ends of $P$.

We say that $G[V(P)\cup \{v_1\ll v_5\}]$ is a {\em jewel} in $G$ 
if $v_1\ll v_5$ are distinct vertices,
 $v_1v_2, v_2v_3,\allowbreak v_3v_4,\allowbreak v_4v_5, v_5v_1$ are edges, $v_1v_3,v_2v_4,v_1v_4$
are nonedges, and $P$ is a path of $G$ between $v_1,v_4$ such that $v_2,v_3,v_5$ have no neighbours
in the interior of $P$. 
(We do not specify whether $v_5$ is adjacent to $v_2,v_3$, but if it is adjacent to one and not the other, then $G$ also contains
a pyramid.)
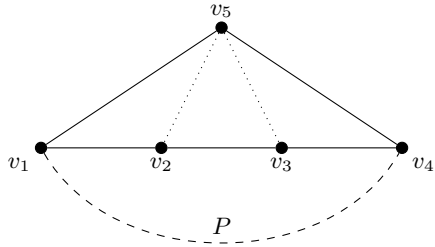
\begin{figure}[H]
\centering

\begin{tikzpicture}[scale=0.8,auto=left]
\tikzstyle{every node}=[inner sep=1.5pt, fill=black,circle,draw]

\node (v1) at (-3,0) {};
\node (v2) at (-1,0) {};
\node (v3) at (1,0) {};
\node (v4) at (3,0) {};
\node (v5) at (0,2) {};
\tikzstyle{every node}=[]
\draw (v1) node [below left]           {\footnotesize$v_1$};
\draw (v2) node [below]           {\footnotesize$v_2$};
\draw (v3) node [below]           {\footnotesize$v_3$};
\draw (v4) node [below right]           {\footnotesize$v_4$};
\draw (v5) node [above]           {\footnotesize$v_5$};
\node  (P) at (0,-1.3)            {\footnotesize$P$};

\draw (v1) -- (v2) -- (v3)--(v4)--(v5)--(v1);
\draw[dotted] (v5) -- (v2);
\draw[dotted] (v5) -- (v3);
\draw[dashed] (v1) to [bend right = 60] (v4);

\end{tikzpicture}

\caption{A jewel. Throughout, dotted lines represent possible edges.} \label{fig:jewel}
\end{figure}

Again, if $G$ contains a jewel then it has an odd hole; and 
it is shown in theorem 3.1 of~\cite{bergealg} that:

\begin{thm}\label{testjewel}
There is an algorithm with the following specifications:
\begin{description}
\item [Input:] A graph $G$.
\item [Output:] Determines whether there is a jewel in $G$.
\item [Running time:] $O(|G|^6)$.
\end{description}
\end{thm}

It is proved in theorem 4.2 of~\cite{bergealg} that:
\begin{thm}\label{testclean}
There is an algorithm with the following specifications:
\begin{description}
\item [Input:] A graph $G$ containing no pyramid or jewel.
\item [Output:] Determines one of the following:
\begin{enumerate}
\item  $G$ contains an odd hole;
\item  there is no clean shortest odd hole in $G$.
\end{enumerate}
\item [Running time:] $O(|G|^4)$.
\end{description}
\end{thm}

Let us say a shortest odd hole $C$ is {\em heavy-cleanable} if there is an edge $uv$ of $C$ such that
every $C$-major vertex is adjacent to one of $u,v$. We deduce:

\begin{thm}\label{testcleanable}
There is an algorithm with the following specifications:
\begin{description}
\item [Input:] A graph $G$ containing no pyramid or jewel.
\item [Output:] Determines one of the following:
\begin{enumerate}
\item  $G$ contains an odd hole;
\item  there is no heavy-cleanable shortest odd hole in $G$.
\end{enumerate}
\item [Running time:] $O(|G|^8)$.
\end{description}
\end{thm}
\Proof
List all the four-vertex induced paths $c_1\dd c_2\dd c_3\dd c_4$ of $G$. For each one, let $X$ be the set of all vertices of $G$
different from $c_1\ll c_4$ and adjacent to one of $c_2,c_3$. We test whether $G\setminus X$ has a clean shortest odd hole, by
\ref{testclean}. If this never succeeds, output that $G$ has no heavy-cleanable shortest odd hole.

To see correctness, note that if $G$ has a
heavy-cleanable shortest odd hole $C$ then $C$ is clean in $G\setminus X$ for some $X$ that we will test (assuming we have
not already detected an odd hole); and when we do so,
\ref{testclean} will detect an odd hole. If $G$ does not have a heavy-cleanable shortest odd hole, then two things might happen:
either \ref{testclean} detects an odd hole for some choice of $X$, or it never detects one. In either case the output is correct.
This proves \ref{testcleanable}.~\bbox

Let us say that a graph $G$ is a {\em candidate} if it has no jewel or pyramid, 
and no heavy-cleanable shortest odd hole (and consequently no hole of length five).
By combining the previous results we deduce:

\begin{thm}\label{testcandidate}
There is an algorithm with the following specifications:
\begin{description}
\item [Input:] A graph $G$.
\item [Output:] Determines one of the following:
\begin{enumerate}
\item  $G$ contains an odd hole;
\item  $G$ is a candidate.
\end{enumerate}
\item [Running time:] $O(|G|^{9})$.
\end{description}
\end{thm}

In view of this, we just need to find a poly-time algorithm to test candidates for odd holes.

\section{Heavy edges}

Let $C$ be a graph that is a cycle, and let $A \subseteq V(C)$.
An {\em $A$-gap} is a subgraph of $C$ composed of a component $X$
of $C\setminus A$, the vertices of $A$ with neighbours in $X$, and the edges between $A$ and $X$.
(So if some two vertices in $A$ are nonadjacent, the $A$-gaps are the paths of $C$ of length $\ge 2$, with both ends
in $A$ and no internal vertex in $A$.) The {\em length} of an $A$-gap is the number of edges in it (so if
$A$ consists just of two adjacent vertices, the $A$-gap has length $|E(C)| - 1$).
We say that $A$ is {\em normal} in $C$ if every $A$-gap is even (and consequently if $C$ has odd length then $A\ne \emptyset$).

The following is proved in theorem 7.6 of~\cite{bergealg}:
\begin{thm}\label{stablenbr}
Let $G$ be a graph containing no jewel or pyramid, let $C$ be a shortest odd hole in $G$, and let
$X$ be a stable set of $C$-major vertices. Then the set of $X$-complete vertices in $C$ is normal.
\end{thm}

In particular, we have:
\begin{thm}\label{manynbrs}
If $G$ is a candidate and $C$ is a shortest odd hole in $G$,
then every $C$-major
vertex has at least four neighbours in $V(C)$.
\end{thm}
\Proof
Let $v$ be $C$-major, and suppose it has at most three neighbours in $V(C)$. Let $A$ be the set of neighbours of $v$ in $C$.
Every $A$-gap is even (since adding $v$ gives a hole 
shorter than $C$), and since $C$ is odd, some edge of $C$ is not in a $A$-gap, that is, 
some two neighbours of $v$ in $V(C)$ are adjacent. Since 
$v$ is $C$-major, it has exactly three neighbours in $V(C)$, and they are not all three consecutive; but then $G[V(C)\cup \{v\}]$ is a pyramid,
a contradiction.
This proves~\ref{manynbrs}.~\bbox

\begin{thm}\label{majorjump}
Let $G$ be a candidate, let $C$ be a shortest odd hole in $G$, and let $x,y$ be nonadjacent $C$-major vertices. Then every induced path between $x,y$ with interior in $V(C)$
has even length.
\end{thm}
\Proof
Let the vertices of $C$ in order be $c_1,c_2\ll c_n,c_1$. By \ref{stablenbr} applied to $\{x,y\}$, some vertex of $C$ is adjacent to both
$x,y$, say $c_n$. Suppose that there is an odd induced path $P$ between $x,y$ with interior in $V(C)$; then $c_n\notin V(P)$, and
since
adding $c_n$ does not give an odd hole (because such an odd hole would be shorter than $C$), it follows that $c_n$
has a neighbour in the interior of $P$. Thus we may assume that the interior of $P$ equals $\{c_1,c_2\ll c_k\}$ for some even $k\ge 2$.
If both $x,y$ have a neighbour in the set $\{c_{k+2}\ll c_{n-2}\}$,
there is an induced path $Q$ between $x,y$
with interior in $\{c_{k+2}\ll c_{n-2}\}$, and its union with one of $x\dd c_n\dd y, P$ is an odd hole of length less than that of $C$,
a contradiction.
Thus one of $x,y$ has no neighbours in $\{c_{k+2}\ll c_{n-2}\}$, say $x$. By \ref{manynbrs}, $x$ has at least four neighbours in $V(C)$;
so it has exactly four, and is adjacent to both $c_{k+1},c_{n-1}$. Hence the neighbours of $x$ in $V(C)$ are $c_{n-1}, c_n, c_{k+1}$, 
and exactly one of $c_1,c_k$. But $k$ is even, since $P$ has odd length, so the path $c_{k+1}\dd c_{k+2}\cc c_{n-1}$ of $C$
is odd; and since adding $x$ does not make an odd hole shorter than $C$, it follows that $c_{k+1}, c_{n-1}$ are
adjacent, and so $k=n-3$. But then the
four neighbours of $x$ in $C$ are consecutive, and so the subgraph induced on $V(C)\cup \{x\}$ is a jewel, a contradiction.
This proves \ref{majorjump}.~\bbox

If $X\subseteq V(G)$, we say an edge $uv$ of $G$ is {\em $X$-heavy} if $u,v\notin X$, and every vertex of $X$ is adjacent to at least
one of $u,v$. We need:

\begin{thm}\label{heavyedge}
Let $G$ be a graph containing no jewel or pyramid or $5$-hole, and 
let $C$ be a shortest odd hole in $G$. Let $X$ be a set of $C$-major vertices, and let $x_0\in X$ be nonadjacent to all other 
members of $X$. Then there is an $X$-heavy edge in $C$.
\end{thm}
\Proof
We proceed by induction on $|X|$. If $X$ is stable then by \ref{stablenbr} some vertex of $C$ is $X$-complete (because the null set is 
not normal), and both edges of $C$ incident with it are $X$-heavy, as required. We assume then that $x_1,x_2\in X$ are adjacent.
From the inductive hypothesis, some edge $u_iv_i$ of $C$
is $(X\setminus \{x_i\})$-heavy, for $i = 1,2$; so we may assume that $x_1$ has no neighbour in $u_2v_2$
(because otherwise $u_2v_2$ is $X$-heavy), and similarly $x_2$ has no neighbour in $u_1v_1$. Consequently $u_1v_1$
and $u_2v_2$ are distinct edges. Since $x_0,x_1$ both have neighbours in $\{u_2,v_2\}$, \ref{majorjump} implies that
they have a common neighbour in $\{u_2,v_2\}$; so we may assume that $x_0,x_1$ are both adjacent to $v_2$, and similarly $x_0, x_2$
are both adjacent to $v_1$. Since the subgraph induced on $\{x_0,x_1,x_2,v_1,v_2\}$ is not a 5-hole, and $v_1\ne v_2$, it follows that
$v_1,v_2$ are adjacent. Since $x_1$ has no neighbour in $\{u_2,v_2\}$, it follows that $u_2\ne v_1$, and $u_1\ne v_2$,
so $u_1,v_1,v_2,u_2$ are in order in $C$. We claim that $v_1v_2$ is $X$-heavy; for suppose not. Then there exists $x\in X$
nonadjacent to $v_1,v_2$; and so $x\ne x_0,x_1,x_2$. Since $u_1v_1$ is $(X\setminus \{x_1\})$-heavy, it follows that $x$
is adjacent to $u_1$, and similarly to $u_2$; but then the subgraph induced on $\{x,u_1,v_1,v_2,u_2\}$ is a 5-hole, a contradiction.
This proves \ref{heavyedge}.~\bbox

\section{The odd holes algorithm}

We can now give the algorithm to detect an odd hole. We first present it in as simple a form as we can, but its running time will be
$O(|G|^{12})$. Then we show that with more care we can bring the running time down to $O(|G|^9)$.

Let $C$ be a hole and $x\in V(G)\setminus V(C)$. An {\em $x$-gap}
is  an induced path of $C$ with length at least two, with both ends adjacent to $x$ and with its internal vertices nonadjacent to $x$. Thus
if $P$ is a $x$-gap then $G[V(P)\cup \{v\}]$ is a hole.
We need the following, theorem 4.1 of~\cite{bergealg}:
\begin{thm}\label{shortpath}
Let $G$ be a graph containing no jewel or pyramid, and let $C$ be a clean shortest odd hole in $G$.
Let $u,v\in V(C)$ be distinct and nonadjacent, and let $L_1,L_2$ be the two subpaths of $C$ joining $u,v$,
where $|E(L_1)| < |E(L_2)|$. Then:
\begin{itemize}
\item  $L_1$ is a shortest path in $G$ between $u,v$, and
\item for every shortest path $P$ in $G$ between $u, v$, $P\cup L_2$ is a shortest odd hole in $G$.
\end{itemize}
\end{thm}

Here then is a preliminary version of the algorithm.
We are given an input graph $G$. First we apply the algorithm of \ref{testcleanable}, and we may assume it determines that
$G$ is a candidate.

Next we enumerate all induced four-vertex paths $c_1\dd c_2\dd c_3\dd c_4$ of $G$, and all three-vertex paths $d_1\dd x\dd d_2$
of $G$. For each choice of $c_1\dd c_2\dd c_3\dd c_4$ and $d_1\dd x\dd d_2$, and each vertex $d_3$ of $G$ (thus we are checking all
$8$-tuples $(c_1,c_2,c_3,c_4,d_1,x,d_2,d_3)$), we do the following:

\begin{itemize}
\item Compute the set $X_1$ of all vertices adjacent to both $d_1$ and $d_2$ that are different from $x$, and compute the set $X_2$ of all vertices that are
adjacent to one of $c_2,c_3$ and different from $c_1,c_2,c_3,c_4$. Let $G'$ be the graph obtained from $G$ by deleting $X_1\cup X_2$.
Compute the set $Y$ of all vertices of $G'$ that are different from and nonadjacent to $x$. 

\item If $d_3\notin Y$, move on to the next $8$-tuple. Otherwise,
check that the distances in $G[Y\cup \{d_1,d_2\}]$ between $d_1,d_3$ and between $d_2,d_3$ are finite and equal (and if not, move on to the next $8$-tuple).

\item For each $y\in Y$, compute the distance in $G[Y\cup \{d_1,d_2\}]$ to $d_1$, to $d_2$ and to $d_3$.
For $i = 1,2$, let $F_i$ be the set of all $y\in Y$ with the sum of the distances to $d_i$ and to $d_3$ minimum; that is, the set of interiors of shortest
paths in $G[Y\cup \{d_1,d_2\}]$ between $d_3$ and $d_i$.
Let $X_3$ be the set of all vertices of $G'$ different from $d_1,d_2,d_3,x$ that are not in $F_1\cup F_2$ and have a neighbour in $F_1\cup F_2\cup \{d_3\}$.

\item Use the algorithm of \ref{testclean} to determine either that $G\setminus (X_1\cup X_2\cup X_3\cup \{x\})$ has an odd hole, or that it has
no clean shortest odd hole. If it finds that there is an odd hole, we output this fact and stop. If after examining all choices of
$8$-tuple we have not found that there is an odd hole, we output that there is none, and stop.

\end{itemize}

Let us see that this algorithm works correctly. Certainly, if the input graph has no odd hole then the output is correct; so we may assume
that 
$G$ is a candidate and $C$ is a shortest odd hole in $G$. Since $C$ is not heavy-cleanable, there is a $C$-major vertex $x$ with an 
$x$-gap of 
length at least three; and so
there is one, $x$ say, with an $x$-gap in $C$ of maximum length, at least three. Let this $x$-gap have ends $d_1,d_2$; 
so $d_1,d_2$ are adjacent to $x$, and there is a path $D$ of $C$ between $d_1,d_2$ such that $x$ has no neighbour in its interior.
Since $x$ is $C$-major, the $C$-distance between $d_1,d_2$ is at least three (because the path of $C$ joining $d_1,d_2$
different from $D$ contains all the neighbours of $v$ in $V(C)$).
From the choice of $x$, every other $C$-major vertex is either adjacent to both $d_1,d_2$, or has a neighbour in the
interior of $D$. Since $C$ is a shortest odd hole, it follows that $D$ has even length; let $d_3$ be its middle vertex. 

By \ref{heavyedge} there is an edge $c_2c_3$ of $C$ such that all $C$-major vertices nonadjacent to $x$
are adjacent to one of $c_2,c_3$. Let $c_1\dd c_2\dd c_3\dd c_4$ be a path of $C$. As the algorithm examines in turn each $8$-tuple,
it eventually will examine the $8$-tuple $(c_1,c_2,c_3,c_4,d_1,x,d_2,d_3)$,
and we will show that in this case
the algorithm will determine that there is an odd hole.

Thus, suppose that the algorithm is now examining the ``correct'' $8$-tuple. Let $X_1,X_2, G',Y$ be as in the first bullet above. 
It follows that $X_1,X_2$ are disjoint from $V(C)$
and so $C$ is a shortest odd hole in $G'$. Since $d_3$ is the middle vertex of $D$, the subpaths of $D$ joining $d_3$ to $d_1,d_2$ both have
length less than $|C|/2$, and so, by \ref{shortpath}, these two subpaths are shortest paths joining their ends with interior in $Y$. Moreover, by \ref{shortpath},
for every choice
of a shortest path $L_i$ in $G'[Y\cup \{d_1,d_2\}]$ joining $d_3,d_i$ (for $i=1,2$), no vertex of $L_i\setminus \{d_i\}$ belongs to or has a neighbour 
in $V(C)\setminus V(D)$; and so the set $X_3$ defined in the third bullet above contains no vertex in $V(C)$. But $X_1\cup X_2\cup X_3\cup \{x\}$
contains every $C$-major vertex, and so $C$ is a clean shortest odd hole in $G\setminus (X_1\cup X_2\cup X_3\cup \{x\})$; and hence when we apply
the algorithm of \ref{testclean} to this subgraph, it will determine that it (and hence $G$) has an odd hole. This completes the proof of correctness.
For the running time: there are $|G|^8$ $8$-tuples to check.  For each one, the sequence of steps above takes time $O(|G|^4)$; and so the total running time is
$O(|G|^{12})$.

Now let us do it more carefully, to reduce the running time. In order to explain the method, let us consider more closely a shortest odd hole $C$
in a candidate $G$. As before, there is a $C$-major vertex $x$ with an $x$-gap of
length at least three; and so
there is one, $x$ say, with an $x$-gap in $C$ of maximum length, at least three. Let this $x$-gap have ends $d_1,d_2$;
so $d_1,d_2$ are adjacent to $x$, and there is a path $D$ of $C$ between $d_1,d_2$ such that $x$ has no neighbour in its interior.
By \ref{heavyedge}, there is an edge $c_2c_3$ of $C$ such that both $x$ and all $C$-major vertices nonadjacent to $x$
are adjacent to one of $c_2,c_3$. Since $x$ is adjacent to one of $c_2,c_3$, not both $c_2,c_3$
belong to the interior of $D$; and since $d_1,d_2$ are nonadjacent, we may assume by exchanging $d_1,d_2$ or $c_2,c_3$
if necessary that
$c_2,c_3,d_1,x,d_2$ are all distinct except that possibly $c_2=d_1$. Now there are six possibilities:
\begin{enumerate}
\item $c_2\ne d_1$ (and hence $c_3$ does not belong to the interior of $D$), and $D$ has length less than $|C|/2$;
\item $c_2\ne d_1$ (and hence $c_3$ does not belong to the interior of $D$), and $D$ has length more than $|C|/2$;
\item $c_2=d_1$, and $c_3$ does not belong to the interior of $D$, and $D$ has length less than $|C|/2$; 
\item $c_2=d_1$, and $c_3$ does not belong to the interior of $D$, and $D$ has length more than $|C|/2$;
\item $c_3$ belongs to the interior of $D$ (and hence $c_2=d_1$), and $D$ has length less than $|C|/2$;  and
\item $c_3$ belongs to the interior of $D$ (and hence $c_2=d_1$), and $D$ has length more than $|C|/2$.
\end{enumerate}
Let us say $C$ is of {\em type $i$} if the $i$th bullet above applies, where $1\le i\le 6$. (Thus $C$ may have more than one type.)
To minimize running time, it seems best to design separate algorithms to test for the six types separately. 
We need the following lemma:
\begin{thm}\label{odddistance}
There is an algorithm with the following specifications:
\begin{description}
\item [Input:] A graph $G$, and two disjoint subsets $A,B$ of $V(G)$, and a vertex $h\notin A\cup B$ with no neighbour 
in $A\cup B$. Also for each $v\in A\cup B$, an induced path $R_v$ between $v$ and $h$, containing no vertex in $A\cup B$ except $v$.
\item [Output:] Determines whether there exist $a\in A$ and $b\in B$ such that $R_a\cup R_b$ is an induced path.
\item [Running time:] $O(|G|^3)$.
\end{description}
\end{thm}
\Proof
We construct a new graph $G'$ by adding edges to $G$,  making $v$ adjacent to every vertex in $V(R_v)\setminus \{v,h\}$,
for each $v\in A\cup B$. Then we test whether there exist $a\in A$ and $b\in B$ with distance four in $G'$.
It is easy to see that a pair $a,b$ has distance four in $G'$ if and only if $R_a\cup R_b$ is an induced path in $G$.
This proves \ref{odddistance}.~\bbox

Now to handle the six types of shortest odd hole. (If $x$ is a vertex of $G$, $N[x]$ denotes the set consisting of $x$ and all its neighbours.)
We begin with:
\begin{thm}\label{type1alg}
There is an algorithm with the following specifications:
\begin{description}
\item [Input:] A candidate $G$.
\item [Output:] Determines either that $G$ has an odd hole, or that $G$ has no shortest odd hole of type $1$.
\item [Running time:] $O(|G|^{9})$.
\end{description}
\end{thm}
\Proof
List all 6-tuples $(c_2,c_3,d_1,x,d_2,d_3)$ of distinct vertices of $G$ such that $c_2c_3$ is an edge, and $d_1\dd x\dd d_2$ is an induced path.
Now we test each such $6$-tuple in turn, as follows.
Let $X_1$ be the set of common  neighbours of $d_1,d_2$ different from $x$, and let $X_2$ be the set of all vertices different from 
$x,c_2,c_3,d_1,d_2$
that are adjacent to one of $c_2,c_3$. Let $C_1$ be the set of all vertices of $G$
different from $c_2,c_3,x$
that are adjacent to $c_2$ and not to $c_3$, and let $C_4$ be the set that are adjacent to $c_3$ and not to $c_2$.
Let $G'$ be the graph obtained from $G$ by deleting $X_1\cup X_2\cup (N[x]\setminus \{d_1,d_2\})$. 

Find the distance in $G'$ between $d_1,d_2$, say $t$. 
If $t$ is infinite
move on to the next $6$-tuple. If $t$ is finite, for each $v\in V(G')$ compute the distance between $v$ and $d_i$ for $i = 1,2$ (setting
the distance to be infinite if there is no path), and let $Y$ be the set of $v\in V(G)$ different from $d_1,d_2$ with the sum of these two distances
equal to $t$. Let $X_3$ be the set of vertices of $G$ different from $x,d_1,d_2$ that are not in $Y$ and have a neighbour in $Y$. 

Let
$G''=G\setminus (X_1\cup X_2\cup X_3\cup \{x\})$. 
For each $v\in C_1\cup C_4$, if there is a path of $G$ between $v$ and $d_3$ such that all its vertices except $v$ belong to $V(G'')$, find such a path
$R_v$ of minimum length.
Let $C_1'$ be the set of $v\in C_1$ such that $R_v$ exists and has even length, and let $C_1''$ be the set of
$v\in C_1$ such that $R_v$ exists and has odd length. Define $C_4', C_4''$ similarly. Apply \ref{odddistance} to test whether there exist $c_1\in C_1'$ and $c_4\in C_4'$ such that the paths $R_{c_1}, R_{c_4}$
are both defined and have union an induced path between $c_1,c_4$ with interior in $V(G'')$. If so output that $G$ has an odd hole. Otherwise apply
\ref{odddistance} again to test whether there exist $c_1\in C_1''$ and $c_4\in C_4''$ such that the paths $R_{c_1}, R_{c_4}$
are both defined and have union an induced path between $c_1,c_4$ with interior in $V(G'')$. If so output that $G$ has an odd hole. 
Otherwise move on to the next $6$-tuple. When all $6$-tuples have been tested, if no odd hole is found, return that $G$ has no shortest odd hole of type 1.

To see the correctness, suppose that $C$ is a shortest odd hole of type 1 in $G$, and let $c_2,c_3,d_1,d_2,x, D$ be as in the definition of type.
Let $d_3$ be the vertex of $C$ that is the middle vertex of the even path $C\setminus \{c_2,c_3\}$. 
When the algorithm tests the $6$-tuple
$(c_2,c_3,d_1,x,d_2,d_3)$, let $X_1,X_2,C_1,C_4$ be as in the description
of the algorithm. From the choice of $c_2,c_3$, every $C$-major vertex nonadjacent to $x$ belongs to $X_2$. 
From the first assertion of \ref{shortpath}, $D$ is a shortest path between $d_1,d_2$ in $G'$, of length $t$ say. By the second assertion of 
\ref{shortpath}, if $L$ is a shortest path between $d_1,d_2$ in $G'$, then no vertex in $V(C)\setminus V(D)$
has a neighbour in $V(L)\setminus \{d_1,d_2\}$; and so $X_3\cap V(C)=\emptyset$. Let $c_1\in C_1$ and $c_4\in C_4$ such that 
$c_1\dd c_2\dd c_3\dd c_4$ is a path of $C$. Thus $C\setminus \{c_1,c_4\}$ is a subgraph of $G''$. Moreover, the paths $R_{c_1},R_{c_4}$
exist, and by the first assertion of \ref{shortpath} they both have length $(|V(C)|-3)/2$; and by the second assertion of \ref{shortpath},
the union of $R_{c_1},R_{c_4}$ is an induced path between $c_1,c_4$. Then adding $c_2,c_3$ gives an odd hole, as required. This proves \ref{type1alg}.~\bbox

The algorithms for type 3 and type 5 are small modifications of this. 

\begin{thm}\label{type3alg}
There is an algorithm with the following specifications:
\begin{description}
\item [Input:] A candidate $G$.
\item [Output:] Determines either that $G$ has an odd hole, or that $G$ has no shortest odd hole of type $3$.
\item [Running time:] $O(|G|^{9})$.
\end{description}
\end{thm}
\Proof
Enumerate all $7$-tuples $(c_1,c_3,c_4,d_1,x,d_2, d_3)$ of distinct vertices such that $c_1\dd d_1\dd c_3\dd c_4$ is an induced path and
$d_1\dd x\dd d_2$ is an induced path. Set $c_2=d_1$. 
For each such $7$-tuple, let $X_1$ be the set of common  neighbours of $d_1,d_2$, and let $X_2$ be the set of vertices different from $c_1,c_2,c_3,c_4$
that are adjacent to one of $c_2,c_3$.
Let $G'$ be the graph obtained from $G$ by deleting $X_1\cup X_2\cup (N[x]\setminus \{c_1,d_2\})$.

Find the distance in $G'$ between $c_1,d_2$, say $t$.
If $t$ is infinite
move on to the next $6$-tuple. If $t$ is finite, for each $v\in V(G')$ compute the distance between $v,c_1$ and between $v,d_2$,
and let $Y$ be the set of $v\in V(G)$ different from $c_1,d_2$ with the sum of these two distances
equal to $t$. Let $X_3$ be the set of vertices of $G$ different from $x,c_1,d_2$ that are not in $Y$ and have a neighbour in $Y$.

Let
$G''=G\setminus (X_1\cup X_2\cup X_3\cup \{x\})$.
Find a shortest path between $c_1,d_3$ in $G''$, and a shortest path between $c_4,d_3$ in $G''$, and test whether their union is an odd induced path
between $c_1,c_4$. If so output that $G$ has an odd hole and otherwise move on to the next $7$-tuple.
If no odd hole is found after testing all $7$-tuples, output that $G$ has no odd hole of type 3.
The proof of correctness is like that for \ref{type1alg} and we omit it. This proves \ref{type3alg}.~\bbox

Similarly we have:
\begin{thm}\label{type5alg}
There is an algorithm with the following specifications:
\begin{description}
\item [Input:] A candidate $G$.
\item [Output:] Determines either that $G$ has an odd hole, or that $G$ has no shortest odd hole of type $5$.
\item [Running time:] $O(|G|^{9})$.
\end{description}
\end{thm}

\begin{thm}\label{type2alg}
There is an algorithm with the following specifications:
\begin{description}
\item [Input:] A candidate $G$.
\item [Output:] Determines either that $G$ has an odd hole, or that $G$ has no shortest odd hole of type $2$.
\item [Running time:] $O(|G|^{9})$.
\end{description}
\end{thm}
\Proof
Enumerate all $6$-tuples $(c_2,c_3,d_1,x,d_2,d_3)$ of distinct vertices where $c_2c_3$ is an edge and $d_1\dd x\dd d_2$ is an induced path.
We test each $6$-tuple in turn as follows.
Let $X_1$ be the set of common  neighbours of $d_1,d_2$ different from $x$, and let $X_2$ be the set of all vertices different from
$x,c_2,c_3,d_1,d_2$
that are adjacent to one of $c_2,c_3$. Let $C_1$ be the set of all vertices of $G$
different from $c_2,c_3,x$
that are adjacent to $c_2$ and not to $c_3,d_3$, and let $C_4$ be the set that are adjacent to $c_3$ and not to $c_2,d_3$.
Let $G'$ be the graph obtained from $G$ by deleting $X_1\cup X_2\cup (N[x]\setminus \{d_1,d_2\})$.
If $d_3\notin V(G')$ move on to the next $6$-tuple.
Find the distance in $G'$ between $d_1,d_3$, and between $d_2,d_3$, and check that they are finite and equal (to $t$ say); if not move on to the next $6$-tuple.
For each $v\in V(G')$ find its distance in $G'$ to $d_1,d_2$ and $d_3$; let $Y_1$ be the set of $v\ne d_1,d_2$ with the sum of its distances to $d_1,d_3$
equal to $t$, and let $Y_2$ be the set of $v\ne d_1,d_2$ with the sum of its distances to $d_2,d_3$ equal to $t$.  (Thus $d_3\in Y_1,Y_2$.)
Let $X_3$ be the set of vertices of $G$ different from $x,d_1,d_2$ that are not in $Y$ and have a neighbour in $Y$.

Let
$G''=G\setminus (X_1\cup X_2\cup X_3\cup \{x\})$.
For each $v\in C_1$, if there is a path of $G$ between $v$ and $d_1$ such that all its vertices except $v$ belong to $V(G'')$, find such a path
$R_v$ of minimum length.
For each $v\in C_4$, if there is a path of $G$ between $v$ and $d_2$ such that all its vertices except $v$ belong to $V(G'')$, find such a path
$R_v$ of minimum length.
Let $C_1'$ be the set of $v\in C_1$ such that $R_v$ exists and has even length, and let $C_1''$ be the set of
$v\in C_1$ such that $R_v$ exists and has odd length. Define $C_4', C_4''$ similarly. Apply \ref{odddistance} in the graph obtained from 
$G[C_1\cup C_4\cup  V(G'')]$ by identifying $d_1,d_2$,  to test whether there exist $c_1\in C_1'$ and $c_4\in C_4'$ such that the paths $R_{c_1}, R_{c_4}$
are both defined and are disjoint and have no edges joining them. If so output that $G$ has an odd hole. Otherwise apply
\ref{odddistance} again to test whether there exist $c_1\in C_1''$ and $c_4\in C_4''$ such that the paths $R_{c_1}, R_{c_4}$
are both defined and are disjoint and have no edges joining them. If so output that $G$ has an odd hole.
Otherwise move on to the next $6$-tuple. When all $6$-tuples have been tested, if no odd hole is found, return that $G$ has no shortest odd hole of type 2.

To see the correctness, suppose that $C$ is a shortest odd hole of type 2 in $G$, and let $c_2,c_3,d_1,d_2,x, D$ be as in the definition of type.
Let $d_3$ be the vertex of $C$ that is the middle vertex of the even path $D$.
When the algorithm tests the $6$-tuple
$(c_2,c_3,d_1,x,d_2,d_3)$, let $X_1,X_2,C_1,C_4$ be as in the description
of the algorithm; then it follows from \ref{shortpath} as before that the output is correct. This proves \ref{type2alg}.~\bbox

Similar modifications handle the remaining two cases, and we omit them. In summary we have:
\begin{thm}\label{finalalg}
There is an algorithm with the following specifications:
\begin{description}
\item [Input:] A candidate $G$.
\item [Output:] Determines whether $G$ has an odd hole.
\item [Running time:] $O(|G|^{9})$.
\end{description}
\end{thm}
\Proof
If $G$ has an odd hole, then it has a shortest odd hole of one of the six types; and by running the algorithms just described we can detect it.~\bbox

Our main result \ref{mainthm} follows immediately from this and \ref{testcandidate}. One final remark: we have an algorithm to 
determine whether $G$ has an odd hole, but what about actually finding an odd hole?  One could obviously do this with an extra 
factor of $|G|$ in the running time, just by deleting vertices and running \ref{finalalg} repeatedly, to find a maximal subset
of the vertex set whose deletion does not destroy all odd holes. But we can do better, and in fact it is easy to adapt the current
algorithm to find an odd hole instead of just detecting the existence of one, with running time $O(|G|^{9})$. We omit the details.

\end{document}